\documentclass[a4paper,12pt]{article}
\usepackage{amsmath}
\usepackage{amssymb}
\usepackage{amscd}
\usepackage{amsthm}
\usepackage{graphicx}
\newtheorem{theorem}{Theorem}
\newtheorem{lemma}[theorem]{Lemma}
\newtheorem*{rtheorem}{Theorem}
\newtheorem{corollary}[theorem]{Corollary}

\newtheorem{remark}[theorem]{Remark}
\usepackage{pdfsync}
\def \PG{\mathrm{PG}}
\def\AG{\mathrm{AG}}
\def\V{\mathrm{V}}
\def\Co{\mathrm{C}}

\def\D{\mathcal{D}}
\def\N{\mathcal{N}}
\def\A{\mathcal{A}}
\def\F{\mathbb{F}}

\def\S{\mathcal{S}}
\def\E{\mathcal{E}}

\def\PGammaL{\mathrm{P}\Gamma\mathrm{L}}
\def\PGL{\mathrm{PGL}}

\title{On sets without tangents \\
and exterior sets of a conic}
\author{Geertrui Van de Voorde}
\date{}
\begin{document}

\maketitle

\begin{abstract} 
A set without tangents in $\PG(2,q)$ is a set of points $\S$ such that no line meets $\S$ in exactly one point. An exterior set of a conic $\mathcal{C}$ is a set of points $\E$ such that all secant lines of $\E$ are external lines of $\mathcal{C}$. 
In this paper, we first recall some known examples of sets without tangents and describe them in terms of determined directions of an affine pointset. We show that the smallest sets without tangents in $\PG(2,5)$ are (up to projective equivalence) of two different types. We generalise the non-trivial type by giving an explicit construction of a set without tangents in $\PG(2,q)$, $q=p^h$, $p>2$ prime, of size $q(q-1)/2-r(q+1)/2$, for all $0\leq r\leq (q-5)/2$. After that, a different description of the same set in $\PG(2,5)$, using exterior sets of a conic, is given and we investigate in which ways a set of exterior points on an external line $L$ of a conic in $\PG(2,q)$ can be extended with an extra point $Q$ to a larger exterior set of $\mathcal{C}$. It turns out that if $q=3$ mod 4, $Q$ has to lie on $L$, whereas if $q=1$ mod $4$, there is a unique point $Q$ not on $L$. 

\end{abstract}
\section{Introduction and preliminaries}
Throughout this paper, $q=p^h$, where $p$ is a prime, and the finite field of order $q$ is denoted by $\F_q$. It is well known that $\F_q$, $q$ odd, contains $(q+1)/2$ squares and $(q-1)/2$ non-squares. The projective plane of order $q$ is denoted by $\PG(2,q)$. If $\S$ is a set of points in $\PG(2,q)$, then a line meeting $\S$ in exactly $1$ point is called a {\em tangent} line (of $\S$), lines containing no points of $\S$ are called {\em external} (lines of $\S$) and lines containing at least $2$ points are called {\em secant} (lines of $\S$). If $\S$ is a set such that no line of $\PG(2,q)$ is a tangent of $\S$, $\S$ is called a {\em set without tangents}.

Let $\mathcal{C}$ be an irreducible conic in $\PG(2,q)$, $q$ odd, then $\mathcal{C}$ has $q+1$ points and a line of $\PG(2,q)$ meets $\mathcal{C}$ in $0$, $1$, or $2$ points. The conic $\mathcal{C}$ has $q+1$ tangent lines, $q(q-1)/2$ external lines and $q(q+1)/2$ secant lines. 
The points, not on $\mathcal{C}$, that lie on a tangent line of $\mathcal{C}$, are called {\em exterior} (points of $\mathcal{C}$) and every exterior point lies on exactly two tangent lines of $\mathcal{C}$. The points not lying on a tangent line of $\mathcal{C}$, are called {\em interior} (points of $\mathcal{C}$). A {\em Desargues configuration} $\mathcal{D}$ in $\PG(2,q)$ is a set consisting of 10 points of $\PG(2,q)$,  and $10$ lines of $\PG(2,q)$, such that every line contains $3$ points of $\mathcal{D}$ and through every point of $\mathcal{D}$ there are 3 lines.

\subsection{Exterior sets of conics}
An {\em exterior set} $\mathcal{E}$ of $\mathcal{C}$ is a set of points such that every secant line of $\mathcal{E}$ is an external line of $\mathcal{C}$.


If $\S$ is a set of $(q+1)/2$ exterior points forming an exterior set of the conic $\mathcal{C}$, then we have the following theorem by Blokhuis, Seress and Wilbrink.
\begin{theorem}{\rm \cite{blokhuis}} Let $\S$ be a set of $(q+1)/2$ exterior points forming an exterior set of the conic $\mathcal{C}$. If $q=1 \mod 4$, $\S$ consists of the $(q+1)/2$ exterior points on an external line of $\mathcal{C}$. If $q= 3 \mod 4$, there exist other examples (at least for $q=7,11,\ldots,31$).
\end{theorem}
It is conjectured by the authors of the same paper (and checked by computer for $q<131$), that only for $q=7,11,\ldots,31$, there exist other examples.

\subsection{Sets without tangents}\label{u}

A non-empty set without tangents in $\PG(2,q)$ is easily seen to have at least $q+2$ points. If $q$ is even, this bound is sharp by the existence of hyperovals. If $q$ is odd, no hyperovals exist and the determination of size of the smallest non-empty set without tangents in $\PG(2,q)$ remains an open problem, except for a few small values of $q$. Let $u_q$ be the size of the smallest set without tangents in $\PG(2,q)$. In \cite{blokhuis2}, Blokhuis, Seress and Wilbrink determined that $u_3=6$, $u_5=10$, $u_7=12$. We showed by computer that $u_9=15$, and $u_{11}=18$.

It is easy to construct a set without tangents of size $2q$: let $\S$ be the set of points on two different lines $L_1$ and $L_2$, different from the intersection point $L_1\cap L_2$. We call this example the {\em trivial} set without tangents.


\subsection{Directions determined by a pointset}

Let $\S$ be a set of points with coordinates $\langle(x,y,z)\rangle$ in $\PG(2,q)$, let $L_\infty$ be the line with equation $z=0$ of $\PG(2,q)$ and let $\PG(2,q)\setminus L_\infty$ be the affine plane $\AG(2,q)$, obtained by removing the line $L_\infty$. The pointset of $\A=\S\setminus L_\infty$ in $\AG(2,q)$ is called the {\em affine part} of $\S$ and consists of points with coordinates $\langle(x_i,y_i,1)\rangle$, $i=1,\ldots, \vert \S\vert$. The set of determined directions $\mathcal{D}$ is defined as
$$\mathcal{D}=\left\{\frac{y_i-y_j}{x_i-x_j}\vert 1\leq i\neq j \leq \vert \S\vert\right\}.$$
We identify a direction $d$ in $\mathcal{D}$ with the point $\langle(1,d,0)\rangle$ on $L$. If $x_i=x_j$ for some $i\neq j$, then the determined direction is $\infty$, which is identified with the point $\langle(0,1,0)\rangle$ on $L_\infty$.

\subsection{Finite geometry codes and stopping sets}

An {\em $\F_p$-linear code} $\Co$ of length $n$ is a vector subspace of $\V(n,p)$. The vectors of $\Co$ are called {\em codewords}. The {\em incidence matrix} of a projective plane $\Pi$ of order $q=p^h$ is the matrix $A$ where the columns are indexed by the points $P_1,\ldots,P_{q^2+q+1}$ of $\Pi$, the rows are indexed by the lines $L_1,\ldots,L_{q^2+q+1}$ of $\Pi$, and with entry $A_{ji}=1$ if $P_j$ lies on $L_i$, and entry $0$ otherwise. The {\em code of $\PG(2,q)$}, denoted by $\Co(2,q)$ is the $\F_p$-span of the rows of $A$. The {\em dual code} $\Co^\bot$ of a code $\Co$ of length $n$ is the vector space consisting of all vectors $v$ of $\V(n,p)$ such that $v.c=0$ for all $c\in \Co$. The {\em support} of a codeword is the set of coordinate positions for which the corresponding entry in the codeword is non-zero. 

The minimum weight of $\Co(2,q)^\bot$ is only known in the case that $q$ is even (then it is $q+2$), and the case that $q$ is a prime $p$ (then it is $2p$). It is not too hard to see that the set of points $\S$ defined by the support of a codeword of $\Co(2,q)^\bot$ is a set without tangents: a codeword has scalar product with all rows of $A$ equal to zero. Since the rows of $A$ correspond to the lines of $\Pi$, this implies that there cannot be lines containing exactly one point of $\S$.
It should be noted that every codeword of $\Co(2,q)^\bot$ gives rise to a set without tangents, but the vice versa part is generally not true. As we have seen, the minimum weight of $\Co(2,p)^\bot$, $p$ prime, is $2p$, but we will give an example of a set without tangents of weight $2p-2$ for $p>5$. However, the smallest known example of a non-trivial set without tangents in $\PG(2,q)$, $q$ not a prime, arises from the support of a codeword in $\Co(2,q)^\bot$.

The dual code of a projective plane is often considered as a so-called {\em LDPC-code} (see \cite{LDPC}), and the number of errors that can be decoded by using iterative decoding over a binary erasure channel is entirely defined by the size of the smallest {\em stopping set} (see \cite{di}). In the case of the LDPC-code of a projective plane $\PG(2,q)$, these {\em stopping sets} are exactly the sets without tangents in $\PG(2,q)$, regardless whether $q$ is even or odd. Hence, the problem of finding the smallest set without tangents in projective planes of odd order is of importance when studying the LDPC-codes of projective planes. This point of view on sets without tangents was used in the papers \cite{stopping2}, \cite{stopping}, where elementary bounds on the size of the smallest stopping sets are derived. However, the results are much weaker than the ones on sets without tangents obtained by Blokhuis, Seress and Wilbrink \cite{blokhuis2}.

\subsection{Arcs, conics, and $\PGammaL(3,q)$}

A {\em $k$-arc} in $\PG(2,q)$ is a set of $k$ points in $\PG(2,q)$ such that no three of them are collinear. It is easy to see that a $k$-arc in $\PG(2,q)$, has at most $q+2$ points. If a $(q+2)$-arc (i.e. a  hyperoval) exists, necessarily $q$ is even. The pointset of an irreducible conic forms a $(q+1)$-arc in $\PG(2,q)$, and if $q$ is odd, Segre showed that the converse also holds.

\begin{theorem}\label{segre} {\sc \cite{segre}}
A $(q+1)$-arc in $\PG(2,q)$, $q$ odd, is an irreducible conic.
\end{theorem}

The following well-known lemmas will be used in the proof of Theorem \ref{hoofd}. The collineation group of $\PG(2,q)$ is denoted by $\PGammaL(3,q)$ and consists of all semi-linear transformations.

\begin{lemma} \label{trans}(1) {\rm (See e.g. \cite[Theorem 2.36]{hughes})} The group $\PGL(3,q)$ acts transitively on the non-empty irreducible conics of $\PG(2,q)$.

(2) {\rm (See e.g. \cite[Theorem 22.6.6]{hirschfeld})} Let $\mathcal{C}$ be an irreducible conic in $\PG(2,q)$, $q$ odd. The stabiliser of $\mathcal{C}$ in $\PGammaL(3,q)$, acts transitively on the external lines of $\mathcal{C}$.
\end{lemma}
\begin{remark} Recall that an irreducible conic in $\PG(2,q)$, $q$ odd, is always non-empty.
\end{remark}

\section{Sets without tangents in $\PG(2,q)$}

\subsection{A lower bound and some old examples}
 If $q$ is odd, the following theorem of Blokhuis, Seress and Wilbrink gives a lower bound on the size of a set without tangents in $\PG(2,q)$.

\begin{theorem}{\rm \cite{blokhuis2} } \label{bl}A non-empty set without tangents in $\PG(2,q)$, $q$ odd, has at least $q+\frac{1}{4}\sqrt{2q}+2$ points.
\end{theorem}
Unfortunately, the bound of Theorem \ref{bl} is probably not sharp and the known examples of sets without tangents have size substantially larger than this lower bound. As we have already seen, the trivial set without tangents in $\PG(2,q)$ contains $2q$ points. In \cite{blokhuis2}, the authors present another example of a set without tangents in $\PG(2,q)$, $q>5$. It has size $2(q-1)$ and arises from two conics: let $\mathcal{C}_1$ be the conic with equation $Z^2=XY$ and let $\mathcal{C}_2$ be the conic with equation $Z^2=aXY$, with $a$ in $\F_q$ such that $1-a$ and $a(a-1)$ are both squares. The points that lie on $\mathcal{C}_1$ or $\mathcal{C}_2$, but not on both, can be shown to be a set without tangents.
For $q$ prime, this example is the best known. If $q$ is not prime, the following construction giving rise to a codeword in $\Co(2,q)^\bot$ by Lavrauw, Storme and Van de Voorde (see also \cite{key2}) improves this bound and is up to our knowledge the smallest known.

\begin{theorem}{\rm \cite{LSV1}} Let $\A$ be the set of points of the form $\{\langle(1,x,x^p)\rangle\vert x\in \F_q\}$ and let $\mathcal{N}$ be the set of points on the line $x=0$ that are not of the form $\langle(0,x,x^p)\rangle$, then $\S\cup \mathcal{N}$ is a set without tangents of size $q+(q-p)/(p-1)$.
\end{theorem}
The set $\A$ in the previous theorem is a set of $q$ affine points and the set $\mathcal{N}$ is the set of non-determined directions of $\A$. This example is the example of smallest size of the following, more general construction.

\begin{theorem} \label{redei}Let $\A$ be a set of $q$ affine points in $\PG(2,q)$, $p>2$, and let $\mathcal{D}$ be the set of determined directions of $\A$, lying on $L_\infty$.  If $\vert \D\vert <(q+3)/2$, then $\A$, together with the complement of $\D$ in $L_\infty$, is a set without tangents. Vice versa, if $\S$ is a set without tangents of size $q+k$ and suppose that there exists a line $L_\infty$ with $k$ points of $\S$. Then the $k$ points on $L_\infty$ are the non-determined directions of the set $\A=\S \setminus L_\infty$.

\end{theorem}
\begin{proof}
If $\vert \D \vert<(q+3)/2$, then the main theorem of \cite{blokhuis1995} implies that every line meets $\A\cup \D$ in $1$ mod $p$ points. Denote the complement of $\D$ in $L_\infty$ by $\N$ and let $N$ be a point of $\N$. Since $N$ is a non-determined direction, and $\vert \A\vert =q$, every line through $N$ meets $\A$ in exactly one point. Let $P$ be a point of $\A$, and let $L$ be a line through $P$, not through a non-determined direction, then, since $L$ meets $\A\cup \D$ in $1$ mod $p$ points, there are at least $p>2$ points of $\A$ in $L$. Hence, $\A\cup \N$ is a set without tangents.

Suppose now that $\S$ is a set without tangents of size $q+k$ and suppose that there exists a line $L_\infty$ with $k$ points of $\S$. If $Q$ is a point of $L_\infty$, not in $\S$, then every line through $Q$ and a point of $\A=\S\setminus L_\infty$, has to contain another point of $\A$, hence, $Q$ is a direction, determined by the points of $\A$. If $Q$ is a point of $L_\infty$ in $\S$, then all $q$ lines through $Q$ have to contain a point of $\S$. Since $\vert A \vert=q$, every line through $Q$ contains exactly one point of $\A$, and hence $Q$ is a non-determined direction.
\end{proof}
\begin{remark} \label{priem}If $q=p$ prime, the set without tangents constructed in the previous lemma is the trivial set without tangents. It also follows that a set without tangents $\mathcal{S}$ in $\PG(2,p)$ of size $\leq 2p$, that contains a line with $p$ points of $\mathcal{S}$, is the trivial set without tangents. For example, this shows that a set  without tangents $\S$ in $\PG(2,3)$ of size 6 is trivial: let $P$ be a point of $\S$. Since there are $4$ lines through $P$ and $5$ points of $\S$ left, there is a $3$-secant $L$ through $P$. 
\end{remark}
\begin{remark}
If we take $\A=\{\langle(1,x,Tr(x))\rangle\}$, then we obtain the set without tangents of size $2q-q/p$, constructed in \cite{blokhuis1995}.
\end{remark}

\subsection{The prime case}
We can exploit the link with determined directions a little bit further to prove that in $\PG(2,p)$, $p$ prime, a set without tangents $\S$, having a secant with `many' points of $\S$, is trivial. For this, we need the following proposition of Ball, which uses the techniques developed by Blokhuis in \cite{blokhuis3}.

\begin{theorem} {\rm\cite[Corollary 4.4]{simeon}} \label{simeon} Let $\AG(2,p)=\PG(2,p)\setminus L_\infty$, and $p$ prime. Let $\mathcal{A}$ be a set of points of $\AG(2,p)$. If there are at least $\vert \A\vert-(p-1)/2$ and at most $p-1$ points $P$ on $L_\infty$ for which the lines through $P$ are all incident with at least one point of $\mathcal{A}$, then $\mathcal{A}$ contains all the points of a line of $\AG(2,p)$.
\end{theorem}

\begin{theorem}Let $\mathcal{S}$ be a set without tangents in $\PG(2,p)$, $p$ prime, with $\vert \S \vert\leq 2p$. If there is a line $L_\infty$ containing $x$ points of $\mathcal{S}$, and $x\geq \vert \S\vert/2-(p-1)/4$, then $\mathcal{S}$ is trivial. 
\end{theorem}

\begin{proof} Let $L_{\infty}$ be the line containing $x$ points of $\mathcal{S}$, with $x\geq \vert \S\vert/2-(p-1)/4$. For every point $N$ in $\S\cap L_\infty$, the $p$ lines through $N$, different from $L_\infty$, contain a point of $\S$. The affine part $\A=\S\setminus L_\infty$, contains $\vert \S\vert-x$ points. If $x=p+1$, then it easily follows that $\vert \S \vert \geq 2p+1$ which is a contradiction. If $x=p$, then by Remark \ref{priem}, $\mathcal{S}$ is trivial. Hence, suppose that $x\leq p-1$. Since $x\geq \vert \S\vert-x-(p-1)/2$, we may apply Theorem \ref{simeon}, and obtain that the affine part of $\mathcal{S}$ contains an affine line. It follows again from Remark \ref{priem} that $\mathcal{S}$ is trivial.
\end{proof}
\begin{corollary} If $\S$ is a set without tangents in $\PG(2,p)$, and $\vert \S \vert<2p$, then a line has at most $\vert \S \vert/2-(p-5)/4$ points of $\S$.
\end{corollary}

\subsection{A new construction}
\begin{lemma}\label{dualconic}
The set of interior points of an irreducible conic $\mathcal{C}$ in $\PG(2,q)$, $q$ odd, $q\geq 5$ is a set without tangents of size $q(q-1)/2$.
\end{lemma}
\begin{proof} Let $\S$ be the set of interior points of $\mathcal{C}$. The point $P$ of $\S$ lies only on secant lines and external lines to $\mathcal{C}$. As a secant line contains $(q-1)/2$ interior points and an external line contains $(q+1)/2$ interior points, every line through $P$ contains at least $(q-3)/2$ other points of $\S$. Since $q\geq 5$, $\S$ is a set without tangents.
\end{proof}

As we will see in the next subsection, this example, together with the trivial example, is the best possible for $\PG(2,5)$. When $q> 5$, we can find examples of smaller size using the same idea.
\begin{theorem}
There exists a set without tangents $\S$ in $\PG(2,q)$, $q$ odd, of size $q(q-1)/2-r(q+1)/2$, for all $0\leq r\leq (q-5)/2$.
\end{theorem}
\begin{proof}
Let $\mathcal{C}$ be an irreducible conic in $\PG(2,q)$, $q$ odd, and let $Q$ be an exterior point of $\mathcal{C}$. Let $\mathcal{R}$ be a set of $r$ external lines through $Q$ to $\mathcal{C}$, where $0\leq r\leq (q-5)/2$. Let $\S$ be the set of internal points to $\mathcal{C}$, not contained in the lines of $\mathcal{R}$. Then $\vert \S\vert=q(q-1)/2-r(q+1)/2$, and every line through a point of $P$ contains at least $(q-3)/2-(q-5)/2=1$ other point of $\S$.
\end{proof}

\subsection{The case $\PG(2,5)$}
Recall from Subsection \ref{u} that $u_5=10$. In this section, we show that every set without tangents in $\PG(2,5)$ of size 10 is either trivial or arises from the example of Lemma \ref{dualconic}.

\begin{rtheorem} Up to projective equivalence, there exist exactly 2 sets without tangents of size 10 in $\PG(2,5)$:
\begin{itemize}
\item[(i)] The set of points on two lines $L_1$ and $L_2$, different from $L_1\cap L_2$.
\item[(ii)] The points on a Desargues configuration, which is the set of internal points of a conic in $\PG(2,5)$.
\end{itemize}

\end{rtheorem}
\begin{proof}
Let $\S$ be a set without tangents of size $10$ in $\PG(2,5)$. There is no $6$-secant of $\S$, since otherwise the number of points in $\S$ would be at least $11$. If $\S$ has a $5$-secant, then by Remark \ref{priem}, $\S$ is trivial. So now assume that $\S$ has no $5$-secants. Let $x_i$ be the number of lines that meet $\S$ in $i$ points. Note that by definition, $x_1=0$.
We count the number of lines in $\PG(2,5)$, the number of couples $(P,L)$ where $L$ is a line through the point $P$ of $\S$, and the number of triples $(P_1,P_2,L)$, where $L$ is a line through the points $P_1$ and $P_2$ of $\S$. This yields
\begin{eqnarray*}
x_0+x_2+x_3+x_4&=&31\\
2x_2+3x_3+4x_4&=&60\\
2x_2+6x_3+12x_4&=&90.
\end{eqnarray*}
It is easy to check that there are only $2$ solutions $(x_0,x_2,x_3,x_4)$ for this system of equations with $x_i\in \mathbb{N}$; they are $(5,21,2,3)$ and $(6,15,10,0)$. Let us assume that we have $x_0=5,x_2=21,x_3=2$, and that $x_4=3$. Let $L_1,L_2,L_3$ be the $4$-secants of $\S$. If two of these lines, say $L_1$ and $L_2$ meet in a point $P$ of $\S$, then there are $4$ other lines through $P$ that have to contain a point of $\S$, so $\vert \S\vert \geq 11$. Hence, $L_1$ and $L_2$ meet in a point $Q$, not in $\S$. The points $P_1,P_2,P_3,P_4$ on $L_3$ all have to be different from the points of $\S$ in $L_1$ and $L_2$, which forces $\S$ to have at least $12$ points, a contradiction.

We conclude that $x_0=6,x_2=15,x_3=10,x_4=0$. Let $L_1,\ldots,L_6$ be the $6$ lines, skew to $\S$. If three of the lines $L_1,\ldots,L_6$ are concurrent in a point $Q$, then the $10$ points of $\S$ have to lie on the three remaining lines through $Q$. Hence, there is a line through $Q$ with at least $4$ points of $\S$ which is a contradiction. This implies that $L_1,\ldots,L_6$ forms a dual $(q+1)$-arc, and thus, by Theorem \ref{segre}, a dual conic. Since there are $6.7/2=21$ points of $\PG(2,5)$ that lie on one of the lines of $\{L_1,\ldots,L_6\}$, the $10$ points of $\S$ are the $31-21$ points, not on one of those lines. The complement of the points on a dual conic, is clearly the set of interior points of that conic. Lemma \ref{dualconic} shows that the set of interior points of a conic is a set without tangents, and since all conics in $\PG(2,q)$, $q$ odd, are projectively equivalent (see Lemma \ref{trans} (1)), all sets of interior points to a conic are projectively equivalent. It is easy to see that the $10$ points of a Desargues configuration form a set without tangents in $\PG(2,5)$. The statement follows by noticing that this set is not the trivial set without tangents, since no line meets the points of a Desargues configuration in $4$ points.
\end{proof}

\section{Exterior sets in $\PG(2,q)$}

If $q=7$ and $q=11$, there are examples of a set without tangents of size 12 and 18 respectively, obtained in the following way: take the $q+1$ points of an irreducible conic, together with $(q+1)/2$ exterior points, no $3$ on a line, forming an exterior set (see \cite{blokhuis}). This construction was the starting point of the investigation of exterior sets of a conic. In the same paper, Blokhuis, Seress and Wilbrink conjecture that if $q>31$, there are no exterior sets consisting of $(q+1)/2$ non-collinear exterior points in $\PG(2,q)$, and they found by computer that for $11<q\leq 31$, all exterior sets consisting of $(q+1)/2$ exterior points contain a line with at least $3$ points of this set. Hence, the cases $q=7$ and $q=11$ are conjectured to be the only cases for which a conic $\mathcal{C}$ and $(q+1)/2$ exterior points of $\mathcal{C}$ form a set without tangents.

We now give a third point of view on the non-trivial set without tangents of size $10$ in $\PG(2,5)$. Let $\mathcal{C}$ be a conic in $\PG(2,5)$ and let $L$ be an external line of $\mathcal{C}$. Let $P_1,P_2,P_3$ be the exterior points on $L$, and denote the external line, different from $L$, through $P_i$ by $M_i$, $i=1,2,3$. It turns out (as we will see in this section) that $M_1,M_2,M_3$ are concurrent, say in $Q$. From this, it easily follows that the points of $\mathcal{C}$, together with $P_1,P_2,P_3,$  and $Q$ form a set without tangents of size $10$. 

In this section, we check whether we can extend this example to a construction of a non-trivial set without tangents of size $2q$ for other values of $q$. For this, as a first step, we need to find a point $Q$ such that $Q$ together with the exterior points on an external line $L$ is an exterior set. As we will see, if $q=3$ mod $4$, this point $Q$ is always contained in the line $L$, and hence, we can never get a set without tangents in this way. If $q=1$ mod $4$ on the other hand, there is a unique point $Q$ satisfying the required condition, not contained in $L$. It follows that we cannot extend this set with other points to find a set without tangents of size $2q$ if $q>5$.


\begin{theorem} \label{hoofd} Let $\mathcal{C}$ be an irreducible conic in $\PG(2,q)$ and let $L$ be an external line of $\mathcal{C}$. Let $\mathcal{S}$ be the set of $(q+1)/2$ exterior points on $L$. If $q=3$ mod $4$, then each point $Q$ such that $\mathcal{S}\cup \{Q\}$ is an exterior set, lies on $L$. If $q=1$ mod $4$, there is a unique point $Q$, not on $L$, such that $\mathcal{S}\cup \{Q\}$ is an exterior set.
\end{theorem}

\begin{proof} By Lemma \ref{trans} (1), we may fix  $\mathcal{C}$ to be the points $\langle(X,Y,Z)\rangle$ satisfying the equation $Y^2=XZ$. 
A line with equation $Z=aX$ is external if and only if $Y^2=aX^2$ has no non-zero solution, hence, if and only if $a$ is a non-square in $\F_q$. By Lemma \ref{trans} (2), we may fix an external line $L$ to be the line with equation $Z=aX$, where $a$ is a fixed non-square in $\F_q$.

The tangent lines of $\mathcal{C}$ are the line $X=0$, together with the lines with equation $\mu^2X-2\mu Y+Z=0$, where $\mu\in \F_q$. The point $P=\langle(0,1,0)\rangle$ lies on $L$ and on the line $X=0$, hence, it is an exterior point of $\mathcal{C}$. The points on $L$, different from $P$, have the form $\langle(1,\xi,a)\rangle$, for some $\xi\in \F_q$. Such a point is exterior if and only if it lies on one of the tangent lines of $\mathcal{C}$, hence, if and only if the equation $\mu^2-2\xi\mu+a=0$ has a solution, which is the case if and only if $\xi^2-a$ is a square.

The lines through $P$ have equation $Z=\lambda X$, and, as before, we see that this line is an external line of $\mathcal{C}$ if and only if $\lambda$ is a non-square. A point $Q$ on the line $Z=\lambda X$ has coordinates $\langle(1,\alpha,\lambda)\rangle$.

The point $Q$ extends $\S$ to an exterior set if and only if the lines connecting $Q$ with the points of $\S$ are external lines. It is clear that if $Q$ is a point on $L\setminus \S$, we find such a set. Hence, we assume that $Q$ does not lie on $L$, so $\lambda\neq a$.

The line $M$ through $\langle(1,\alpha,\lambda)\rangle$ and $\langle(1,\xi,a)\rangle$ has equation
$$\ (\alpha a-\xi\lambda)X+(\lambda-a)Y+(\xi-\alpha)Z=0.$$  
We see that the line $M$ is external if and only if $(\lambda-a)^2-4(\alpha a-\xi\lambda)(\xi-\alpha)$ is a non-square. 

Suppose there exists a point $Q=\langle(1,\alpha,\lambda)\rangle$, not on $L$  extending the set $\mathcal{S}$ to an exterior set, (hence $\lambda\neq a$ is a non-square). Then, for all $\xi$ with $\xi^2-a$ a square ($\ast$), it holds that $(\lambda-a)^2-4(\xi-\alpha)(\alpha a -\xi\lambda)$ is a non-square $(\ast \ast)$. Since $(\ast)$ is satisfied if and only if the point $\langle(1,\xi,a)\rangle$ is an exterior point on $L$, and there are $(q+1)/2$ exterior points on $L$, and $\langle(0,1,0)\rangle$ is one of them, the number of elements $\xi$ for which $(\ast)$ holds is $(q-1)/2$. The condition $(\ast \ast)$ is satisfied if and only if the line connecting $\langle(1,\alpha,\lambda)\rangle$ with $\langle(1,\xi,a)\rangle$ is an external line. A point $\langle(1,\alpha,\lambda)\rangle$ lies on at most $(q+1)/2$ external lines, of which the line $Z=\lambda X$, through $\langle(0,1,0)\rangle$ is one. Hence, at most $(q-1)/2$ elements of $\F_q$ satisfy condition $(\ast \ast)$, and from our assumption, we get that the $(q-1)/2$ elements of $\F_q$ satisfying condition $(\ast)$ are exactly those for which $(\ast \ast)$ holds.
If $x$ is an element of $\F_q$ for which $(\ast)$ holds, then also $-x$ satisfies $(\ast)$ and if $y$ is an element of $\F_q$ for which $(\ast \ast)$ holds, then $\frac{\alpha(a+\lambda)}{\lambda}-y$ also satisfies $(\ast \ast)$. Since the elements satisfying $(\ast)$ and $(\ast \ast)$ are the same, we get that the (not necessarily distinct) values 

\begin{eqnarray*}
x&,&-x\\
\frac{\alpha(a+\lambda)}{\lambda}-x&,&x-\frac{\alpha(a+\lambda)}{\lambda},\\
\frac{2\alpha(a+\lambda)}{\lambda}-x&,&x-\frac{2\alpha(a+\lambda)}{\lambda},\\
\frac{3\alpha(a+\lambda)}{\lambda}-x&,&x-\frac{3\alpha(a+\lambda)}{\lambda},\ldots,\\
\frac{(p-1)\alpha(a+\lambda)}{\lambda}-x&,&x-\frac{(p-1)\alpha(a+\lambda)}{\lambda}
\end{eqnarray*}

 are satisfying $(\ast)$. 
Suppose that two of the above $2p$ values coincide, then either $x=0$ or $\frac{\alpha(a+\lambda)}{\lambda}=0$. If $x=0$ and $\frac{\alpha(a+\lambda)}{\lambda}\neq 0$, the number of different values is $p$. The above argument holds for every solution $x$ of $(\ast)$. Hence, if $\frac{\alpha(a+\lambda)}{\lambda}\neq 0$, the number of different values for which $(\ast)$ holds is a multiple of $p$, but $(p^h-1)/2$ is not a multiple of $p$ which is a contradiction. We conclude that $\frac{\alpha(a+\lambda)}{\lambda}=0$, so either $\alpha=0$ or $\lambda=-a$.

First assume that $q=3$ mod $4$, then $-1$ is a non-square. This implies that $-a$ is a square, so if $\lambda=-a$, $\lambda$ is a non-square, and we obtain a contradiction. If $\alpha=0$ then the line $Y=0$ through $\langle(1,0,\lambda)\rangle$ and the exterior point $\langle(1,0,a)\rangle$ on $L$ meets $\mathcal{C}$ in $\langle(1,0,0)\rangle$ and $\langle(0,0,1)\rangle$, hence it is not an external line. This proves the theorem in the case that $q=3$ mod $4$.

If $q=1$ mod $4$, then $-1$ is a square. First note that $Q=\langle(1,0,-a)\rangle$ extends $\S$ to a larger exterior set since $\xi^2-a$ is a square if and only if $4a(a-\xi^2)$ is a non-square. We will now show that the point $Q=\langle(1,0,-a)\rangle$ is the unique point satisfying this condition. By the previous argument, we only need to check the points of the form $\langle(1,\alpha,-a)\rangle$ and $\langle(1,0,\lambda)\rangle$.

Suppose that the point $Q=\langle(1,\alpha,-a)\rangle$, with $\alpha\neq 0$ extends $\S$ to a larger exterior set. Then $\xi^2-a$ is a square if and only if $4a(a-\xi^2+\alpha^2)$ is a non-square, hence, if and only if $(\xi^2-a-\alpha^2)$ is a square. It also follows that 
$\xi^2-a$ is a non-square if and only if $(\xi^2-a-\alpha^2)$ is a non-square, hence, $(\xi-a)(\xi^2-a-\alpha^2)$ is always a square. There are $(q+1)/2$ squares in $\F_q$, so there are $(q+1)/2$ values $x$ for which $\xi^2-a=x$. 
We claim that $X(X-\alpha^2)=\nu Z^2$ ($\ast \ast \ast$) with $\nu$ a non-square has $(q+1)/2$ different solutions for $X$. The number of solutions $X$ to
$(\ast \ast \ast)$ is the number of $X$'s for which $\langle(X,1,Z)\rangle$ is a solution of $X^2-\alpha^2 X Y-\nu Z^2=0$. If $\alpha\neq 0$, the last equation is the equation of an irreducible conic $\mathcal{C}'$, and it is clear that none of the $q+1$ points on $\mathcal{C}'$ has  $Y=0$, so we may take $Y=1$. Since a fixed value of $X$ occurs in at most $2$ points, there are at least $(q+1)/2$ different solutions for $X$.
Hence, there is at least one $x$ for which $\xi^2-a=x$, and $x(x-\alpha^2)=\nu z^2$ for some $z$, so $(\xi^2-a)(\xi^2-a-\alpha^2)$ is a non-square, a contradiction. We conclude that $\alpha=0$.

Now suppose that the point $Q=\langle(1,0,\lambda)\rangle$, with $\lambda\neq -a$ extends the set of exterior points on the line $L$ to a larger exterior set, then $\xi^2-a$ is a square if and only if $(\lambda-a)^2+4\lambda \xi^2$ is a non-square and vice versa. Hence, $(\xi^2-a)((\lambda-a)^2+4\lambda \xi^2)$ is a non-square. There are $(q+1)/2$ different values for $X=\xi^2$, and with a similar argument as before, we get that $(X-a)((\lambda-a)^2+4\lambda X)=Z^2$ has $(q+1)/2$ different solutions for $X$, since $(X-aY)((\lambda-a)^2Y+4\lambda X)=Z^2$ is the equation of an irreducible conic if and only if $\lambda\neq-a$. So we conclude that $\lambda=-a$, which implies that the point $Q=\langle(1,0,-a)\rangle$ extending $\S$ to a larger exterior set is unique.
\end{proof}

{\bf Acknowledgment:} The author would like to thank the anonymous referees for their helpful suggestions.


\end{document}